\documentclass[11pt,reqno]{amsart}
\usepackage{amsmath,amssymb}

 \makeatletter
 \evensidemargin  \oddsidemargin
 \makeatother

\begin{document}
\thispagestyle{empty}
 \setcounter{page}{1}

\begin{center}
{\large\bf CONE NORMED SPACES}

\vskip.20in

{\bf M. Eshaghi Gordji,  M. Ramezani, H. Khodaei and H. Baghani} \\[2mm]

\address{Department of Mathematics,
Semnan University\\ P. O. Box 35195-363, Semnan, Iran}

\email{madjid.eshaghi@gmail.com,ramezanimaryam873@gmail.com,
 khodaei.hamid.math@gmail.com, h.baghani@gmail.com}

\vskip 5mm

\end{center}

\vskip 5mm
 \noindent{\footnotesize{\bf Abstract.}
 In this paper, we introduce the cone normed spaces and cone bounded linear
 mappings. Among other things, we prove  the Baire category theorem  and  the
 Banach--Steinhaus theorem  in cone normed spaces.

 \vskip.10in

  \newtheorem{df}{Definition}[section]
  \newtheorem{rk}[df]{Remark}
   \newtheorem{lem}[df]{Lemma}
   \newtheorem{thm}[df]{Theorem}
   \newtheorem{pro}[df]{Proposition}
   \newtheorem{cor}[df]{Corollary}
   \newtheorem{ex}[df]{Example}

 \setcounter{section}{0}
 \numberwithin{equation}{section}

\vskip .2in

\begin{center}
\section{Introduction}
\end{center}
Let $E$ be a Banach space and $P$ be a subset of $E$. $P$ is
called a cone whenever

 (1) $P$ is a closed, non-empty set and $P\not =\{0\}$,

 (2) $ax+by\in P$ for all $x,y\in P$ and $a,b\geq 0$,

 (3) $P\cap(-P)=\{0\}$.\\
   For a given cone $P\subseteq E$, we can define a partial
   ordering $\leq$ with respect to $P$ by $x\leq y$ if and only if
   $y-x\in P$. We shall write $x<y$ if $x\leq y$ and $x\not= y$,
   while $x\ll y$ will stands for $y-x\in int P$, where $int P$
   denoted the interior of $P$. The cone $P$ is normal if there is
   a number $M>0$ such that for all $x,y\in E$
   $$0\leq x\leq y \hspace{0.5 cm} \Longrightarrow\hspace{0.5 cm}  \|x\|\leq M\|y\|.$$
   The least positive number satisfying the above is called the
   normal constant of $P$ \cite{Hu}. It is clear that $ M\geq 1$.
     In the following we always suppose that  $E$ is a real Banach
     space and $P$ is a cone in $E$ with
    $int P\not=\varnothing$ and $\leq$ is a partial ordering with
     respect to $P$.

\begin{df} (\cite{Hu})
Let $X$ be non-empty set . Suppose that the mapping $d: X\times X\to
E$ satisfies:

 (1)$~0\leq d(x,y)$ for all $x,y\in X$ and
$d(x,y)=0$ if and only if $x=y$,

(2) $d(x,y)=d(y,x)$ for all $x,y\in X$,

 (3) $ d(x,y)\leq
d(x,z)+d(z,y)$ for all $x,y,z\in X.$
\end{df}
Then $d$ is called a cone metric on $X$, and $(X,d)$ is called
cone metric space.\\
\begin{ex}
Let $E={\Bbb R}^2$, $P=\{(x,y)\in E : x,y\geq 0\}$, $X=\Bbb R$ and
$d:X \times X \to E$ defined by $d(x,y)=(|{x-y}|,\alpha|x-y|)$,
where $\alpha \geq 0$ is a constant. Then $P$ is a normal cone
with normal constant $M=1$ and $(X,d)$ is a cone metric
space\cite{Re}.
\end{ex}
\begin{df} (\cite{Hu}) Let $(X,d)$ be a cone metric space, $x\in X$ and
$\{x_n\}$ a sequence in $X$. Then

 (1) $\{x_n\}$ is said to
convergent to $x$ whenever for every $c\in E$, with $0\ll c$ there
is a positive integer $N$ such that $d(x_n,x)\ll c$ for all $n\geq
N$. We denote this by $\lim_{n\to\infty} x_n=x$ or $x_n\to x$ as
$n\to \infty$.

 (2) $\{x_n\}$ is said to be a Cauchy sequence
whenever for every $c\in E$ with $0\ll c$ there is a positive
integer $N$ such that $d(x_n,x_m)\ll c$ for all $n,m\geq N$.

 (3) $(X,d)$ is called a complete cone metric space if every Cauchy
sequence is convergent.
\end{df}
 Let us recall [1] that if $P$ is a normal cone, then
$\{x_n\}\subseteq X$ converges to $x$ if and only if $d(x_n,x)\to
0$ as $n\to\infty$. Furthermore, $\{x_n\}\subseteq X$ is a Cauchy
sequence if and only if $d(x_n,x_m)\to 0$ as $n,m\to\infty$.
\vskip5mm
\begin{df} Let $(X,d)$ be a cone metric space and $B\subseteq X.$

(1) A point $b\in B$ is called an interior point of $B$ whenever
there exists a point $p$, $0\ll p$, such that $B_p(b)\subseteq B$
where $B_p(b):=\{y\in X : d(b,y)\ll p\}$.

 (2) A subset $B\subseteq
X$ is called open if each element of $B$
is an interior point of $B$.\\
The family  $\beta=\{B_e(x) : x\in X, 0\ll e\}$ is a sub-basis for
a topology on $X$. We denote this cone topology by $\tau_c$. The
topology $\tau_c$ is a Hausdorff and first countable \cite{Re}.
\end{df}
In this paper we suppose that $P$ is a normal cone with normal
constant $M$ and fixed $c_0$ with $0\ll c_0$.
\begin{center}
\section{ Cone Normed spaces}
\end{center}
\begin{df}Let $X$ be real vector space. Suppose that the mapping
$\|.\|_p :X\to E$ satisfies:

 (1)$~\|x\|_p\geq 0$ for all $x\in X$
and $\|x\|_p=0$ if and only if $x=0$ ,

 (2) $\|\alpha x\|_p=|\alpha|~\|x\|_p$, for all $x\in X$ and
$\alpha\in\Bbb R~$,

 (3)$\|x+y\|_p\leq \|x\|_p+\|y\|_p$ for all
$x,y\in X$. \\Then $\|.\|_p$ is called a cone norm on $X$ and
$(X,\|.\|_p)$ is called a cone normed space.
\end{df}

It is easy to see that every normed space is a cone normed space by
putting $E:=\mathcal R$, $P:=[0.\infty)$.

\begin{ex} Let $E=l_1$, $P=\{\{x_n\}\in E : x_n\geq 0, for~
all ~n\}$ and $(X,\|.\|)$ be a normed space and $\|.\|_p :X\to E$
defined by $\|x\|_p={\{\frac{\|x\|}{2^n}\}}$. Then $P$ is a normal
cone with constant normal $M=1$ and $(X,\|.\|_p)$ is a cone normed
space.\\\
 Let $(X,\|.\|_p)$ be a cone normed space. Set
$d(x,y)=\|x-y\|_p$. It is easy to see that $(X,d)$ is a cone
metric space. $d$ is called  getting cone metric of cone norm
$\|.\|_p$.
\end{ex}
\begin{df}We say that the cone normed space $(X,\|.\|_p)$ is a
cone Banach space when getting cone  metric of  $\|.\|_p$ is
complete.
\end{df}
\begin{lem} Let $(X, \|.\|_p)$ and $(Y,\|.\|_p)$ be cone normed
spaces, and let  $T$ be a linear map from $X$ into $Y$. If  $~T$
has any one of the five following properties, it has all five of
them:\\
(a) (continuity at a point) For some fixed $x_0\in X$ we have:
Given $0\ll c$ there is a $0\ll t$ such that $\|Tx-Tx_0\|_p\ll c$
whenever $\|x-x_0\|_p\ll t$.\\
(b) (continuity at zero) For $0\ll c$ there is a $0\ll t$ such
that $\|Tx\|_p\ll c$ whenever $x\in X$ and $\|x\|_p\ll t$.\\
(c) (continuity at every point of x) For any $x\in X$ we have:
Given $0\ll c$ there is a $0\ll t$ such that $\|Tx-Ty\|_p\ll c$
whenever $y\in X$ and $\|y-x\|_p\ll t$.\\
(d) (uniform continuity) Given $0\ll c$ there is a $0\ll t$ such
that $\|Tx-Ty\|_p\ll c$ whenever $x,y\in X$ and $\|x-y\|_p\ll
t$.\\
(e) (sequential continuity) Given any sequence $\{x_n\}\subseteq
X$ which is convergent to a point $x_0\in X$, the sequence
$\{Tx_n\}\subseteq Y$ is convergent to the point $Tx_0\in Y$.
\end{lem}
\paragraph{\bf Proof:} First assume that $T$ has property (a). So
for some $x_0\in X$  and any $0\ll c$ we can choose $0\ll t$ such
that $\|Tx-Tx_0\|_p\ll c$ whenever $\|x-x_0\|_p\ll t$. Then for
any $w\in X$ with $\|w\|_p\ll t$ we have $\|T(w+x_0)-Tx_0\|\ll c$
because $\|(w+x_0)-x_0\|_p\ll t$. But $T$ is  linear this says
that $\|Tw\|_p\ll c$ whenever $\|w\|_p\ll t$ and we have shown
that (a) implies (b).\\
~~ Now suppose that $T$ has property (b), let $x\in X$ and $0\ll
c$ be given. There is a $0\ll t$ such that $\|Tw\|_p\ll c$
whenever $\|w\|_p\ll t$. We have $\|T(y-x)\|_p\ll c$ whenever
$\|y-x\|_p\ll t$; just use $y-x$ in place  of $w$. Again recalling
that $T$ is linear we see that (b) implies (c). Clearly (c)
implies (a). Thus (a), (b) and (c) are equivalent.\\
~~Let us show that (b) implies (d). Given $0\ll c$ we may choose
$0\ll t$ so that $\|w\|_p\ll t$ implies $\|Tw\|_p\ll c$. Now if
$x,y\in X$ and $\|x-y\|_p \ll t$ then $\|T(x-y)\|_p\ll c$. Since
$T$ is linear (b) implies (d) and clearly (d) implies (b). Thus
(a) through (d) are equivalent.\\
We will complete the proof by showing that (b) and (e) are
equivalent. First suppose that $T$ has property (b) and let $0\ll
c$ be given. Then we have choose $0\ll t$ that $\|Tw\|_p\ll c$
whenever $\|w\|_p\ll t$. Now suppose $\{x_n\}\subseteq X$
converges $x_0$. Then we can find positive integer $N$ such that
$\|x_n-x_0\|_p\ll t$ whenever $n\geq N$. Thus $\|T(x_n-x_0)\|_p\ll
c$ whenever $n\geq N$. Clearly this says that $\{Tx_n\}$ converges
to $Tx_0$. Now assume (e) and negation of (b). So we are supposing
that there is a $0\ll c$ such that for any $0\ll t$ we can find
$w_t\in X$ with $\|w_t\|_p\ll t$ and $\|Tw_t\|_p\not\ll c$. Thus
for this $c$ we can find $\{w_n\}\subseteq X$ such that
$\|w_n\|_p\ll \frac{c}{2^n}$ and $\|Tw_n\|_p\not \ll c$ for all
$n$. But $\{w_n\}$ is converges to $0$ because
$\|~~{\|w_n\|}_p\|\leq \frac{M \|c\|}{2^n}\to 0$ as $n\to\infty$.
By (e) $\{Tw_n\}$ must converge to $T0=0$ and this impossible.
\hfill{$\Box$}

\
\begin{pro}Let $(X,\|.\|_p)$ be a cone normed space, and $x\in
X$, $0\ll c$. Then \\
$$y\in\overline{B_c(x)}\Longleftrightarrow(\exists\{z_n\}\subseteq{B_c(x)}\hspace{0.1 cm};\hspace{0.5 cm}z_n\to
y).$$
\end{pro}
\paragraph{\bf Proof:}Let $y\in\overline{B_c(x)}$. Then for any
positive integer $n$, $z_n\in B_{\frac{c}{2^n}}(y)\cap
B_c(x)\not=\varnothing$. We obtain $z_n\to y$ as $n\to \infty$. Now
we suppose that $\{z_n\}\in B_c(x)$ is a sequence that $z_n\to y$ as
$n\to \infty$. Let $W$ be a open set such that $W$ consists of
 $y$. There is $0\ll p$ such that $B_p(y)\subseteq W$. We choose
the positive integer $n$ such that $\|z_n-y\|_p\ll p$. Hence,
$z_n\in B_p(y)$ and $W\cap B_c(x)\not=\varnothing$. So
$y\in\overline{B_c(x)}$.\hfill$\Box$\\\\
 We need the following lemma to prove the Baire category theorem.
\begin{lem}Let $(X,\|.\|_p)$ be a cone normed space, $x\in X$
and $0\ll c$. Then $\overline{B_\frac c 2(x)}\subseteq B_c(x)$.
\end{lem}
\paragraph{\bf Proof:} Let $y\in\overline{B_\frac c 2(x)}$. Then
there is a sequence $\{z_n\}\subseteq {B_\frac c 2(x)}$ such that
$z_n\to y$ as $n\to \infty$. We can choose the positive integer
$n$ such that $\|z_n-y\|_p\ll\frac c 2$. We obtain that
$\|x-y\|_p\leq \|x-z_n\|_p+\|z_n-y\|_p\ll \frac c 2+\frac c 2= c$.
Letting $a=(\|x-z_n\|_p+\|z_n-y\|_p)-\|x-y\|_p$, by attention that
action $+$ is continuous in $E$, we have
$$c-\|x-y\|_p=(c-(\|x-z_n\|_p+\|z_n-y\|_p)+~a\in a+ int P=int
(a+P)\subseteq int P.$$
 So $\|x-y\|_p\ll c$ and  thus  $y\in
B_c(x)$.\hfill$\Box$
\begin{thm}(Baire category theorem) Let
$(X, \|.\|_p)$ be a cone Banach space. Then every countable
intersection of dense and open sets is dense.
\end{thm}

\paragraph{\bf Proof:} Let $\{A_n\}\subseteq X$ be a sequence of
dense and open sets. Suppose $x\in X$, and $W$ is a open set such
that $W$ consists of  $x$. Then there is $0\ll r$ such that
$B_r(x)\subseteq W$. Since $A_1$ is dense in $X$, we obtain
$z_1\in A_1\cap B_r(x)\not=\varnothing$. But $A_1\cap B_r(x)$ is
open and hence there exists $0\ll r'$ such that
$B_{r'}(z_1)\subseteq A_1\cap B_r(x)$. We can choose the positive
real number $k_1$ such that $k_1<min\{\frac 1 2,\frac 1 {2
\|{r'}_1\|}\}$. By setting $r_1=k_1{r'}_1$, we have
$r_1\ll\frac{r'}{2}$ and $\overline{B_{r_1}(z_1)}\subseteq
B_{r'_1}(z_1)$. Since $A_2$ is dense in $X$, we have $z_2\in
A_2\cap B_{r_1}(z_1)\not=\varnothing$, but this set is open. Thus
there exists $0\ll {r'}_2$ such that $B_{r'_2}(z_2)\subseteq
A_2\cap B_{r_1}(z_1)$. By choosing $0<k_2<min\{\frac 1 2,\frac 1
{2^2 \|{r'}_2\|}\}$, we have $r_2=k_2{r'}_2\ll\frac{{r'}_2} {2}$
and hence $\overline{B_{r_2}(z_2)}\subseteq B_{r'_2}(z_2)$. Since
$A_3$ is dense in $X$ then let $z_3\in A_3\cap B_{r_2}(z_2)$.
Since $A_3$ is open, there is a $0\ll {r'}_3$ such that
$B_{r'_3}(z_3)\subseteq A_3\cap B_{r_2}(z_2)$. We can choose the
real number $0<k_3$ for which $k_3<min\{\frac 1 2,\frac 1 {2^3
\|{r'}_3\|}\}$. Setting $r_3=k_3{r'}_3$, we conclude that
$r_3\ll\frac{r'_3}{2}$ and $\overline
{B_{r_3}(z_3)}\subseteq {B_{r'_3}(z_3)}$. Repeating the above argument we obtain \\
$$...\overline {B_{r_3}(z_3)}\subseteq \overline {B_{r_2}(z_2)}\subseteq
\overline {B_{r_1}(z_1)}.$$
 We claim $r_n\to 0$, because $r_n=k_n r'_n\ll\frac 1 {2^n\|r'_n\|}~ r'_n$.
 So $\|r_n\|\leq \frac M
{2^n\|r'_n\|}\|r'_n\|=\frac M {2^n}\to 0$ as $n\to \infty$.
Moreover , we can show that $\{z_n\}$ is Cauchy sequence in $X$.
To see this, let $\epsilon >0$ is given, there is a positive
integer $N$ such that $M\|r_n\| <\epsilon$ for all $n\geq N$. So
$\| \|z_m-z_n\|_p \|\leq M\|r_n\|<\epsilon$ for all $m>n\geq N$.
This means that $\{z_n\}$ is Cauchy sequence. Since $X$ is cone
Banach space, there is $z\in X$ such that
$\displaystyle\lim_{n\to\infty}z_n= z$. For any positive integer
$N$, if $n>N$ then $z_n\in B_{r_N}(z_N)$, so $z\in\overline
{B_{r_N}(z_N)}$ and hence
 $$z\in \bigcap_{N=1}^{\infty} B_{r_N}(z_N) \subseteq
\bigcap_{N=1}^{\infty}A_N\bigcap B_r(x)\subseteq
\bigcap_{N=1}^{\infty}A_N\bigcap W.$$ This says that
$\bigcap_{N=1}^{\infty}A_N$ is dense in $X$.\hfill$\Box$
\begin{cor}
Every cone Banach space is of second category.
\end{cor}

\begin{df}A subset $A\subseteq E$ is upper bounded, if there
exists $0\leq t$ such that $a\leq t$ for all $a\in A$, $t$ is an
upper bound for $A$. We say that $P$ has supremum property, if for
every upper bounded set $A$ in $P$ least upper bound exists in
$P$, we show this element to $sup A$.
\end{df}
\begin{ex}
Let $E=\Bbb R^2$ , $P=\{(x,y)\in E : x,y\geq 0\}$. $P$ is a normal
cone with normal constant $M=1$, and $P$ has supremum property.
Because  if $A\subseteq P$ is upper bounded set, then $sup
A=(sup_{z\in A} \pi_1(z),sup_{z\in A} \pi_2(z))$. Where $\pi_1$
and $\pi_2$ are projections on first and second components,
respectively.
\end{ex}
From now on, we suppose that $P$ has supremum property.
\begin{df} Let $(X,\|.\|_p)$ be a cone normed space. A subset $A$ in $X$ is cone bounded,
if $\{\|x\|_p ; x\in A\}$ is upper bounded.
\end{df}
\begin{df} Let $(X,\|.\|_p)$ and $(Y,\|.\|_p)$ be cone normed
spaces, and $\Lambda : X\to Y$ be a linear mapping. $\Lambda$ is
cone bounded if the set $\Lambda (B_{c_0}(0)$) is a cone bounded
set.\\
 Let $B(X,Y)$ denotes the set of all cone bounded linear mappings
from $X$ into $Y$. It is easy to see that $\|\Lambda\|_p=sup
\{\|\Lambda(x)\|_p : \|x\|_p\ll c_0\}$ is a cone norm on $B(X,Y)$.
\end{df}
In the  following, we obtain a linear mapping that is not cone
bounded.
\begin{ex} Let $E=\Bbb R^2$, $P=\{(x,y)\in E : x,y\geq 0\}$,
$c_0=(1,1)$ and let $X$ be the set of all real- valued polynomials
on interval $[0,1]$ and $\|.\|_u$ is supremum norm on $X$, that is
$\|f\|_u= sup\{|f(x)|: x\in [0,1]\}$ for all $f\in X$. Let
$\|.\|_p: X\to E$ defined by $\|f\|_p=\|f\|_u c_0$. It is easy to
see that $(X,\|.\|_p)$ is a cone normed space. Suppose $D : X\to
X$ defined by $D(f)=f'$. Then $D$ is a linear mapping that is not
cone bounded.
\end{ex}

\begin{thm}(The Uniform Boundness Principle)
Let $(X,\|.\|_p)$ be a cone Banach space and $(Y,\|.\|_p)$ be a
cone normed space. Suppose $A \subseteq B(X,Y)$ is pointwise
bounded, that is for each $x\in X$, the set $\{Tx : T\in A\}$ is
cone bounded. Then $A$ is a cone bounded set in $B(X,Y)$.
\end{thm}

\paragraph{\bf Proof:}
Let
$$ E_n=\{x\in X :\|Tx\|_p\leq n c_0\ for all~~ T\in A\}$$
for all $n\in \Bbb N$. It is easy too see that the set $\{y\in Y
:\|y\|_p\leq n c_0\}$ is closed. Hence, for each positive integer
$n$, $E_n$ is a closed set. We claim $X=\bigcup_{n=0}^\infty E_n$.
Too see this, since $0\ll c_0$, choose  $0<\delta$ such that
$$c_0+\{x\in E : \|x\|\leq \delta\}\subseteq P.$$
If $x\in X$ is given. Setting $\alpha_x=sup\{\|Tx\|_p : T\in A\}$.
Choose a positive integer $n$ such that
$\|\frac{\alpha_x}{n}\|<\delta$. So $c_0-\frac{\alpha_x}{n}\in
c_0+\{x\in E : \|x\|\leq \delta\}\subseteq P$ and $\alpha_x\ll
nc_0$. This shows that $x\in E_n$. So $X=\bigcup_{n=0}^\infty
E_n$. But $X$ is a cone Banach space  and hance $X$ is  of second
category. Then there is a positive integer $k$ such that
$int{E_k}=int\overline {E_k}\not=\varnothing$. Suppose $y\in int
{E_k}$. We can choose $0\ll r$ such that $B_r(y)\subseteq E_k$.
Repeating the above method, there exists a positive integer $m$
such that $\displaystyle\frac {c_0} m\ll r$. Now, we let $T\in A$
and $x\in X$, $\|x\|_p\ll c_0$. Then $\displaystyle\|(y+\frac x m)
- y\|_p=\|\frac x m \|_p\ll \frac {c_0} m \ll r$. This means that
$y+\frac x m\in B_r(y)$. We have
\begin{align*}
\|Tx\|_p=m \|T\frac x m\|_p&=m\|T(y+\frac x m)- Ty\|_p\\
&\leq m(\|T(y+\frac x m)\|_p+\|Ty\|_p) \\
&\leq m(k+k)=2mk.
\end{align*}
Thus $\|T\|_p\leq 2mk$. In the other words $A$ is a cone bounded
set in $B(X,Y)$.\hfill$\Box$\\\\

The question arises here is  whether  Hahn--Banach theorem and open
mapping theorem can be extended similar to cone normed spaces
 or not?


{\small


\end{document}